\documentclass[12pt]{article}

\renewcommand{\baselinestretch}{1.5}
\usepackage[intlimits]{amsmath} 
\usepackage{amsthm,mathrsfs}    
\usepackage[comma,authoryear]{natbib}
\usepackage{url}
\numberwithin{equation}{section}
\usepackage{amssymb,amsmath}
\usepackage{subfigure}
\usepackage{bm}
\usepackage{booktabs} 
\usepackage{graphicx,graphics,epsfig}
\usepackage{float}

\theoremstyle{plain}
\newtheorem{thm}{Theorem}

\newcommand{\bthm}{\begin{thm}}
\newcommand{\ethm}{\end{thm}}

\newcommand{\bpf}{\begin{proof}}
\newcommand{\epf}{\end{proof}}

\theoremstyle{definition}

\newcommand{\YS}{\operatorname{YS}}

\usepackage[margin=.9in,a4paper]{geometry}
\usepackage[compact]{titlesec}
\newcommand{\bib}{\bibliography{C://Research/Bibtex/ref-bib}\bibliographystyle{Chicago}}
\usepackage{deep-macro}

\begin{document}
\begin{center}
{\Large {\bf Nonlinear Time Series Modeling:\\[.3em] A
Unified Perspective, Algorithm, and Application}}
\\[.2in]
Subhadeep Mukhopadhyay\\
\emph{Department of Statistical Science, Temple University}\\
Emanuel Parzen\\
\emph{Department of Statistics, Texas A\&M University}
\end{center}
\begin{abstract}
A new comprehensive approach to nonlinear time series analysis and modeling is developed in the present paper. We introduce novel data-specific mid-distribution based Legendre Polynomial (LP) \emph{like} nonlinear transformations of the original time series $\{Y(t)\}$ that enables us to adapt all the existing stationary linear Gaussian time series modeling strategy and made it applicable for non-Gaussian and nonlinear processes in a robust fashion. The emphasis of the present paper is on empirical time series modeling via the algorithm LPTime. We demonstrate the effectiveness of our theoretical framework using daily S\&P 500 return data between Jan/2/1963 - Dec/31/2009. Our proposed LPTime algorithm systematically discovers all the `stylized facts' of the financial time series automatically \emph{all at once}, which were previously noted by many researchers one at a time.
\end{abstract}
\vskip1em

\vskip.65em
\renewcommand{\baselinestretch}{.65}
\setlength{\parskip}{.65ex}
{\small
\tableofcontents
}
\renewcommand{\baselinestretch}{1.2}
\setlength{\parskip}{1.25ex}
\section{Introduction}

When one observes a sample $Y(t), t=1, \ldots, T$, of a (discrete parameter) time series $Y(t)$, one seeks to nonparametrically learn from the data a stochastic model with two purposes: (a1) scientific understanding; (a2) forecasting (predict future values of the time series under the assumption that the future obeys the same laws as the past). Our prime focus in this paper is on developing a nonparametric empirical modeling technique for nonlinear (stationary) time series that can be used by data scientists as a practical tool for obtaining insights into (i) the temporal dynamic patterns and (ii) the internal data generating mechanism--a crucial step for achieving (a1) and (a2).
\vskip.4em

Under the assumption that the time series is stationary (which can be extended to asymptotically stationary) the distribution of $Y(t)$ is identical for all $t$, the joint distribution of $Y(t)$ and $Y(t+h)$ depends only on lag $h$. Typical estimation goals are as follows:
\vskip.5em

(1) \emph{Marginal modeling.} The identification of marginal probability law (in particular, the heavy tailed marginal densities) of a time series plays a vital role in financial econometrics. Notations: Common quantile $Q$, inverse of distribution function $F$, respectively denoted $Q(u;Y), 0<u<1$ and $F(y;Y)$. Mid-distribution is defined as $\Fm(y;Y)=F(y;Y)-.5 \Pr(Y(t)=y)$.
\vskip.5em
(2) \emph{Correlation modeling.} Covariance function (defined for positive and negative lag $h$) $R(h;Y)=\Cov[Y(t),Y(t+h)]$. $R(0;Y)=\var[Y(t)]$, $\mu=\Ex[Y(t)]$ assumed $0$ in our prediction theory. Correlation function $\rho(h)=\Cor[Y(t),Y(t+h)]=R(h;Y)/R(0;Y)$.
\vskip.5em
(3) \emph{Frequency-domain modeling.} When covariance is absolutely summable, define spectral density function $f(\om;Y)=\sum R(h;Y) \, e^{- 2 \pi i \om h}, -1/2 < \om < 1/2$.
\vskip.5em
(4) \emph{Time-domain modeling.} Time domain model is linear filter relating $Y(t)$ to white noise $\ep(t)$, $\cN(0,1)$ independent random variables. Autoregressive scheme of order $m$, a predominant linear time series technique for modeling conditional mean is defined as (assuming $\Ex[Y(t)]=0$ )
\beq
Y(t) - a(1;m)Y(t-1) - \ldots -  a(m;m) Y(t-m)\,=\,\sigma_m \ep(t),
\eeq
with the spectral density function given by
\beq \label{eq:arsdf}
f(\om;Y)~=~\dfrac{\sigma^2_m}{\big|   1 - \sum_{k=1}^m a(k;m) e^{2\pi i \om k}    \big|^2}.
\eeq
To fit an AR model, compute linear predictor of $Y(t)$ given $Y(t-j),j=1,\ldots,m$ by
\beq
Y^{\mu,m}[t]\,=\,\Ex \big[ Y(t) \mid Y(t-1),\ldots,Y(t-m) \big] \,=\, a(1;m) Y(t-1) + \cdots +a(m;m)Y(t-m).
\eeq
Verify that the prediction error $Y[t] - Y^{\mu,m}[t]$ are white noise. Best fitting AR order is identified by Akaike criterion AIC (or
Schwarz's criterion BIC) as value of $m$ minimizes
\[\rm{AIC}(m)\,=\,2 \log \si_m + 2 m/n.\]
In what follows we aim to develop a parallel modeling framework for non-linear time series.
\section{From Linear to Nonlinear Modeling} 
Our approach to nonlinear modeling, called LPTime, is via approximate calculation of conditional expectation $\Ex[Y(t)|  Y(t-1),\ldots,Y(t-m)]$. Because with probability one $Q(F(Y)) = Y$, one can prove that the conditional expectation of $Y(t)$ given past values $Y(t-j)$ is equal to (with probability one) conditional expectation of $Y(t)$ given past values $\Fm(Y(t-j))$, which can be approximated by linear orthogonal series expansion in score functions $T_k[\Fm(Y(t-j))]$ constructed by Gram Schmidt orthonormalization of powers of
\beq
T_1 \,=\,\dfrac{\Fm(Y(t);Y) - .5}{\si[\Fm(Y(t);Y)]},
\eeq
where $\si[\Fm(Y(t);Y)]$ si the standard deviation of the mid-distribution transform random variable given by $\sqrt{(1-\sum_y p^3(y))/12}$, $p(y)$ denotes the probability mass function of $Y$. This score polynomials allows us to \emph{simultaneously} tackle the discrete (say count-valued) and continuous time series. Note that for $Y$ is \emph{continuous}, The $T_1$ reduces to:
\beq
T_1\,=\,\sqrt{12} \big( F(Y(t)) - .5 \big).
\eeq
and all the higher-order polynomials $T_j$ can be compactly expressed as $\Leg_j[F(Y)]$, where $\Leg_j(u),0<u<1$ denotes orthonormal Legendre polynomials. It is worthwhile to note that $T_j$ are orthonormal \emph{polynomials of mid-ranks} (instead of polynomials of the original $y$'s), which \emph{injects robustness} into our analysis while allowing to capture nonlinear patterns. Having constructed score functions of $y$ denoted by $T_j$, we transform into unit interval by letting $y=Q(u;Y)$ and defining
\beq
S_j(u;Y)\,=\, T_j\big[\Fm(Q(u;Y))\big],~~T_j(y;Y)\,=\,S_j\big[ \Fm(Y(t))\big].
\eeq
In general our score functions are custom constructed for each distribution function $F$ which can be discrete or continuous.


\begin{figure*}[!htb]
 \centering
 \includegraphics[height=.5\textheight,width=\textwidth,keepaspectratio,trim=1.5cm .5cm 1.5cm 1.5cm]{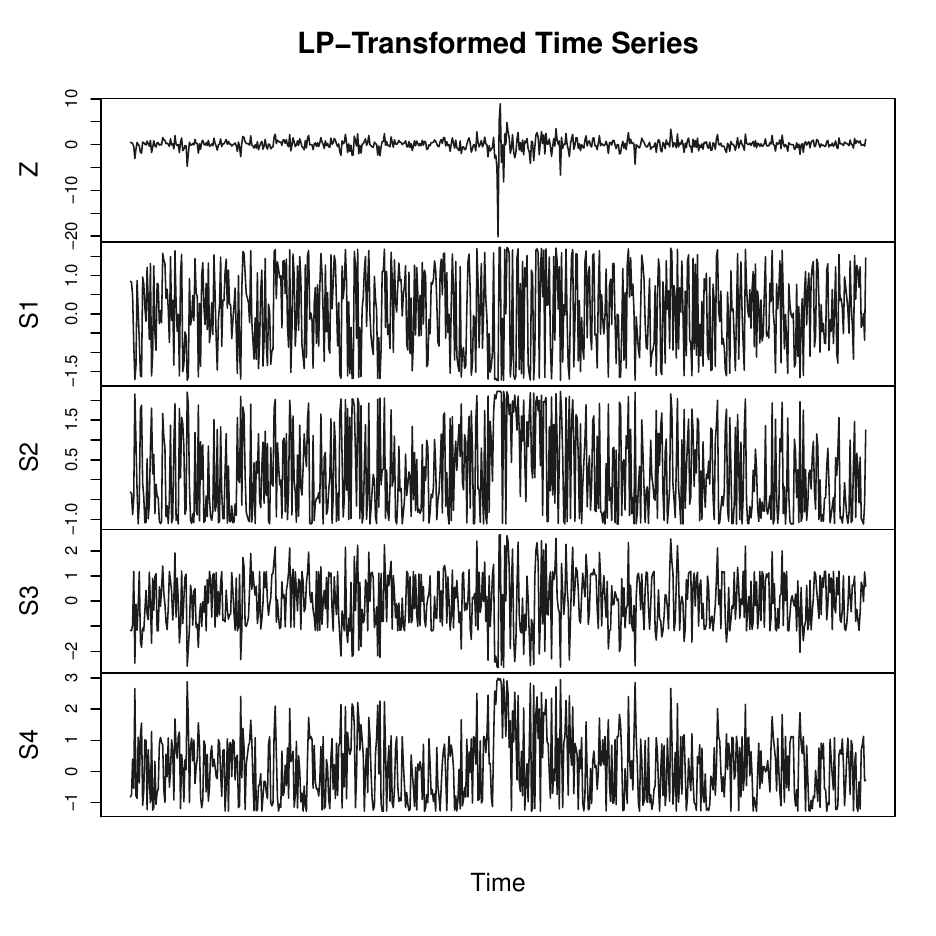} \\
\vspace{-.5em}
\caption{LP-transformed S\&P 500 daily stock returns between Oct 1986 and Oct 1988. This is just a small part of the full time series from Jan/2/1963 - Dec/31/2009 (cf. Sec 3.1).} \label{fig:dis} %
\end{figure*}
\vskip.15em

\section{Nonparametric LPTime Analysis}
Our LPTime empirical time series modeling strategy to nonlinear modeling of a univariate time series $Y(t)$ is based on linear modeling of the multivariate time series:
\beq
\rm{Vec}(\YS)(t)\,=\,\big[ \YS_1(t), \ldots,   \YS_k(t) \big]^{T},
\eeq
where $\YS_k(t)=T_k[\Fm(Y(t))]$, our tailor-made orthonormal mid-rank based nonlinear transformed series.
We summarize below the main steps of the algorithm LPTime. To better understand the functionality and applicability of LPTime we break it into several inter-connected steps each of which highlights:
\begin{description}
  \item[(a)] Algorithmic modeling aspect [how it works]
  \item[(b)] Required theoretical ideas and notions [why it works]
  \item[(c)] Application to daily S\&P 500  return data between Jan/2/1963 - Dec/31/2009 [empirical proof-of-work].
\end{description}
\subsection{The Data and LP-Transformation}
The data used in this paper is daily S\&P 500 return data between Jan/2/1963 - Dec/31/2009 (defined as $\log(P_t/P_{t-1})$ where $P_t$ is the closing price on trading day $t$). We begin our modeling process by transforming the given univariate time series $\{Y(t)\}$ into multiple (robust) time series by means of a special data-analytic construction rules described in (2.1)-(2.3) and (3.1).
We display the original ``normalized'' time series  $\cZ(Y(t))=(Y(t) - \Ex[Y(t)])/\si[Y(t)]$ and the transformed time series $\YS_1(t), \ldots,   \YS_k(t)$ on a single plot.

Fig \ref{fig:dis} shows the first look at the transformed S\&P 500 return data between Oct 1986 and Oct 1988. These newly constructed time series works as a universal preprocessor for any time series modeling in contrast with other adhoc power transformations. In the next sections we will describe how the temporal patterns of these multivariate LP-transformed series $\rm{Vec}(\YS)(t)=\{\YS_1(t), \ldots,   \YS_k(t)\}$  generate various insights for the time series $\{Y(t)\}$ in a organized fashion.
\begin{figure*}[t]
 \centering
 \includegraphics[height=.22\textheight,width=.32\textwidth,trim=1cm .25cm .5cm 1.2cm]{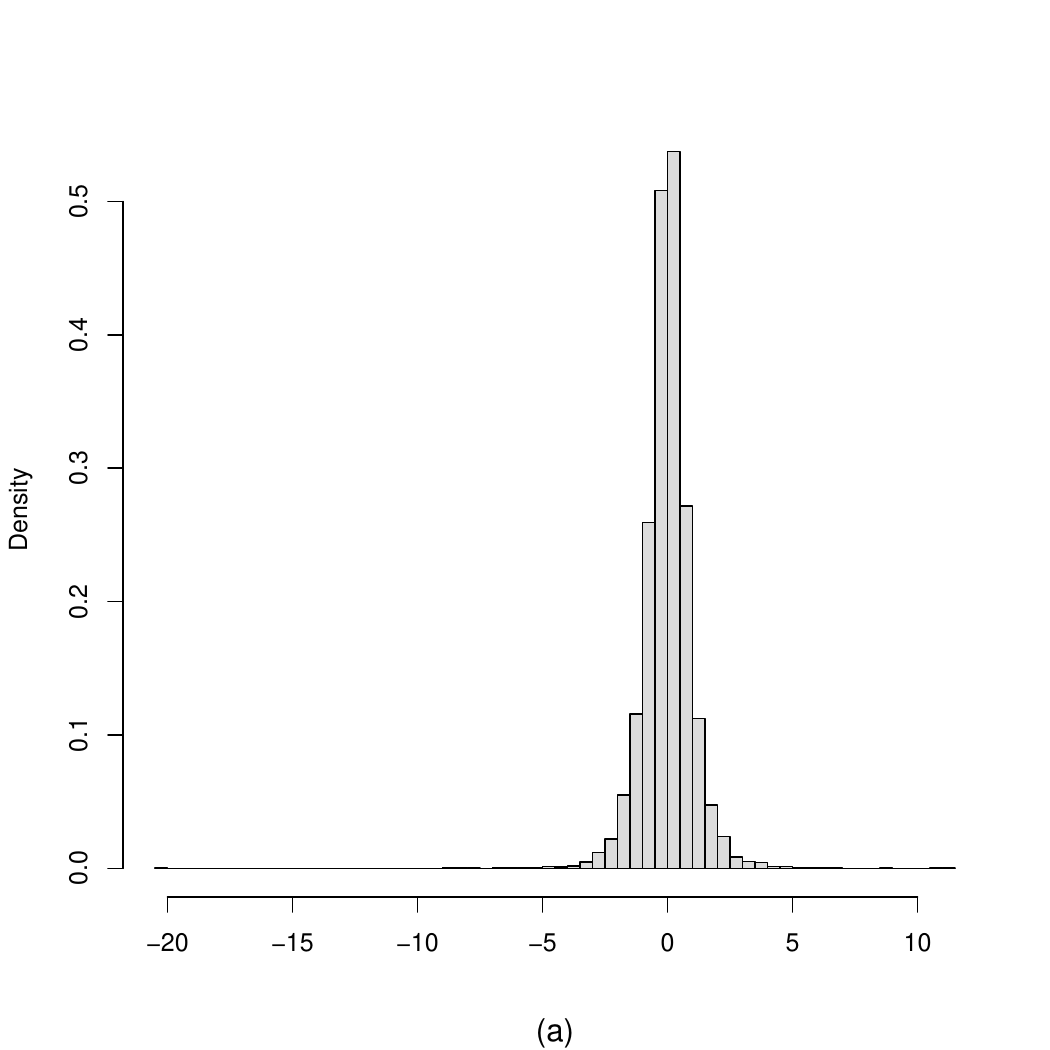}
  \includegraphics[height=.22\textheight,width=.32\textwidth,trim=.75cm .25cm .75cm 1.2cm]{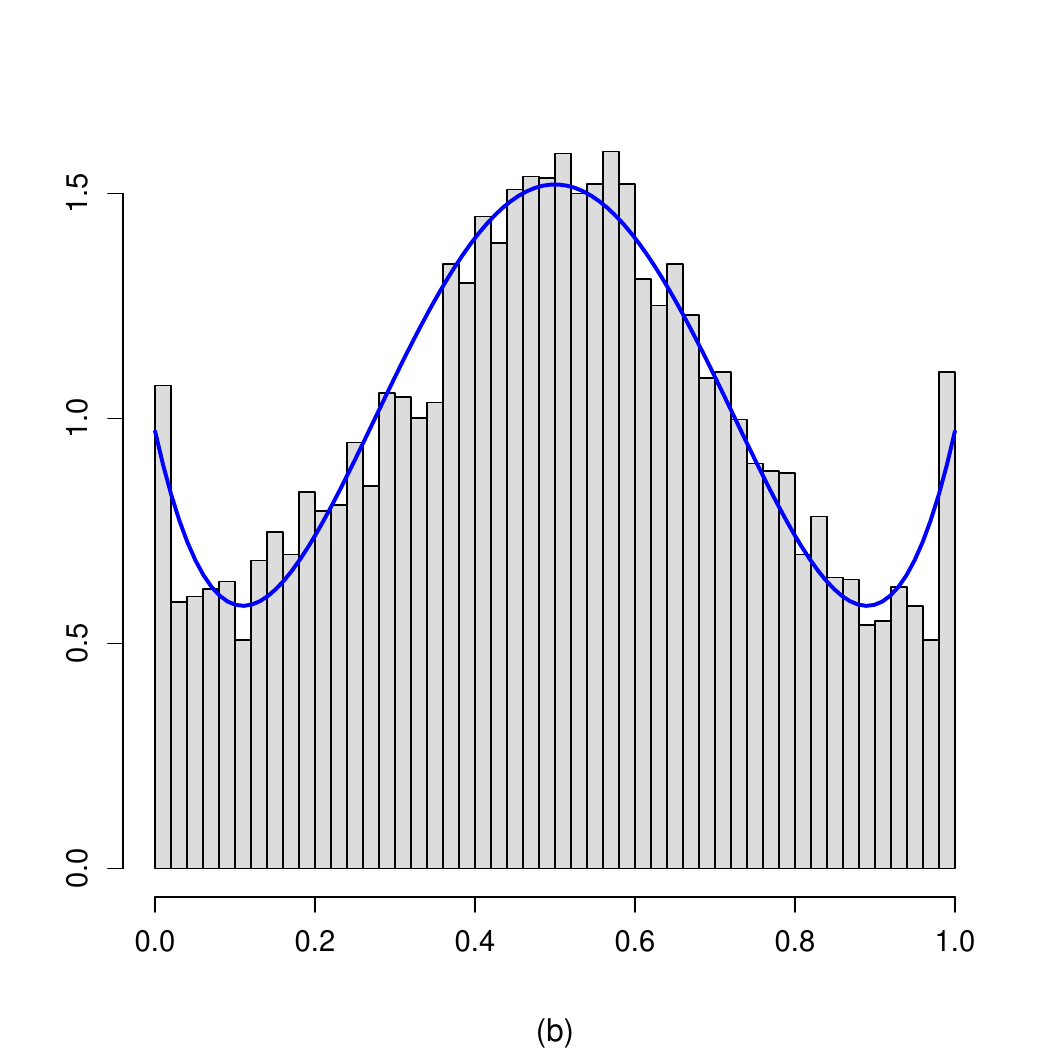}
  \includegraphics[height=.22\textheight,width=.32\textwidth,trim=.5cm .25cm 1cm 1.2cm]{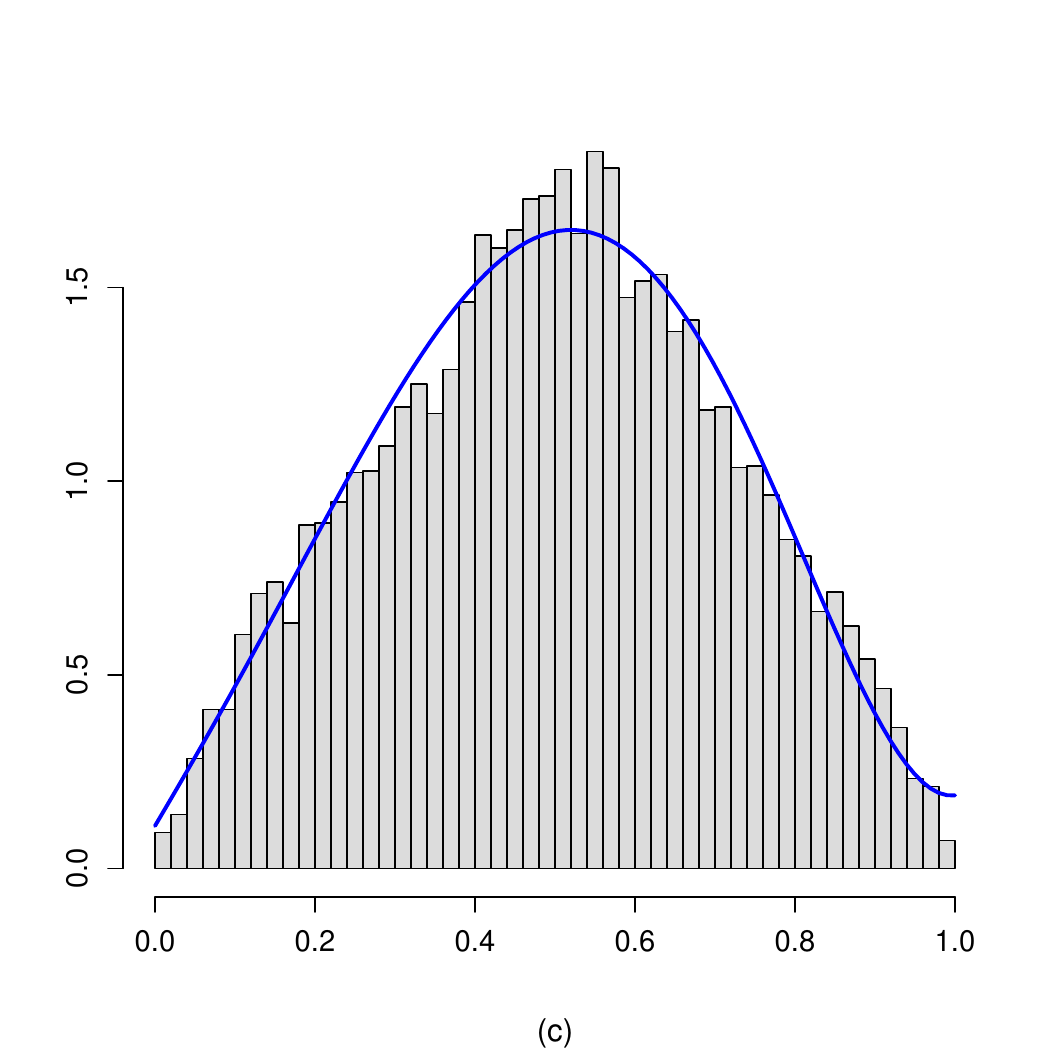}
\vspace{-.5em}
\caption{(a) The marginal distribution of daily returns; (b) plots the histogram of $\Phi(y_i)$ and display the LP-estimated comparison density curve. and (c) shows the associated comparison density estimate with $G$ as $t$-distribution with $2$ degrees of freedom.} \label{fig:LPt-marginal}
\end{figure*}
\subsection{Marginal Modeling}
Our time series modeling starts with the identification of probability distributions.
%
%
\subsubsection{Non-normality Diagnosis}
Does the Normal probability distribution provide a good fit to the S\&P 500 return data? Fig \ref{fig:LPt-marginal}(a) clearly indicates distribution of daily return is certainly non-normal. At this point the natural question is \emph{how} the distribution is different from the assumed normal? A quick insight into this question can be gained by looking at the distribution of the random variable $U=G(Y)$, called comparison density \citep{parzen97,Deep17LPMode} given by:
\beq d(u;G,F)\,=\,\dfrac{f(Q(u;G))}{g(Q(u;G))},~~~~~~0\le u \le1, \eeq
where $Q(u;G)=\inf\{x: G(x) \ge u\}$ is the quantile function. The \emph{flat uniform shape} of the estimated comparison density provides a quick graphical diagnostic to test the fit of the parametric $G$ to the true unknown distribution $F$. The Legendre polynomial based orthogonal series comparison density estimator is given by
\beq d(u;G,F)=1+\sum_j \LP[j;G,F] \Leg_j(u),~~0<u<1\eeq
where the Fourier coefficients $\LP[j;G,F]=\Ex[\Leg_j \circ\, G(Y)]$.

\vskip.3em
For $G=\Phi$, Fig \ref{fig:LPt-marginal}(b) displays the histogram of $U_i=\Phi(Y_i)$ for $i=1,\ldots,n$. The corresponding comparison density estimate $\widehat{d}(u;G,F)\,=\,1-0.271\Leg_2(u) - 0.021 \Leg_3(u) + 0.193\Leg_4(u)$ is shown in blue curve, which reflects the fact that distribution of daily return has \emph{(i) sharp peaked (inverted ``U'' shape),  (ii) negatively skewed with (iii) fatter tails} than the Gaussian distribution. We can carry our similar analysis by asking whether t-distribution with degrees of freedom $2$ provides a better fit. Fig \ref{fig:LPt-marginal}(c) demonstrates the full analysis, where the estimated comparison density $\widehat{d}(u;G,F)\,=\,1-0.492\Leg_2(u) -0.015 \Leg_3(u) +0.084\Leg_4(u)$ indicates (iv) t-distribution fits the data better than normal, specially in the tails, although not a \emph{fully} adequate model.
\vskip.3em
The \emph{shape} of comparison density (along with the histogram of $U_i=G(Y_i)$, $i=1,\ldots, n$) captures and exposes the adequacy of the assumed model $G$ for the true unknown $F$--thus act as an exploratory as well as confirmatory tool.

\subsection{Copula Dependence Modeling}
Distinguishing uncorrelatedness and independence by properly quantifying association is an essential task in empirical nonlinear time series modeling.
\subsubsection{Nonparametric Serial Copula}
We display the nonparametrically estimated smooth serial copula density $\cop(u,v;Y(t),Y(t+h))$ to get much finer understanding of the lagged interdependence structure of a stationary time series. For continuous distribution define the copula density for the pair $(Y(t),Y(t+h))$ as the joint density of $U=F(Y(t))$ and $V=F(Y(t+h))$, which is estimated by sample mid-distribution transform
$\tilde U=\tFm(Y(t))$,\, $\tilde V=\tFm(Y(t+h))$. Following \cite{Deep14LP} and \cite{D12b}, we expand copula density (square integrable) in a orthogonal series of product LP-basis functions as
\beq \label{eq:cop}
\cop(u,v;Y(t),Y(t+h))\,-\,1~=~\sum_{j,k} \LP[j,k;Y(t),Y(t+h)]\,S_j(u;Y(t))\,S_k(v;Y(t+h)),
\eeq
where $S_j(u;Y(t))=\YS_j(Q(u;Y(t));Y(t))$. The Eq \ref{eq:cop} allows us to pictorially represent the information present in the LP-comoment matrix via copula density. The various ``shapes'' of copula density gives insight into the structure and dynamics of the time series.

\begin{figure*}[!htb]
 \centering
 \includegraphics[height=.4\textheight,width=\textwidth,keepaspectratio,trim=2cm 1.2cm .8cm 1.8cm]{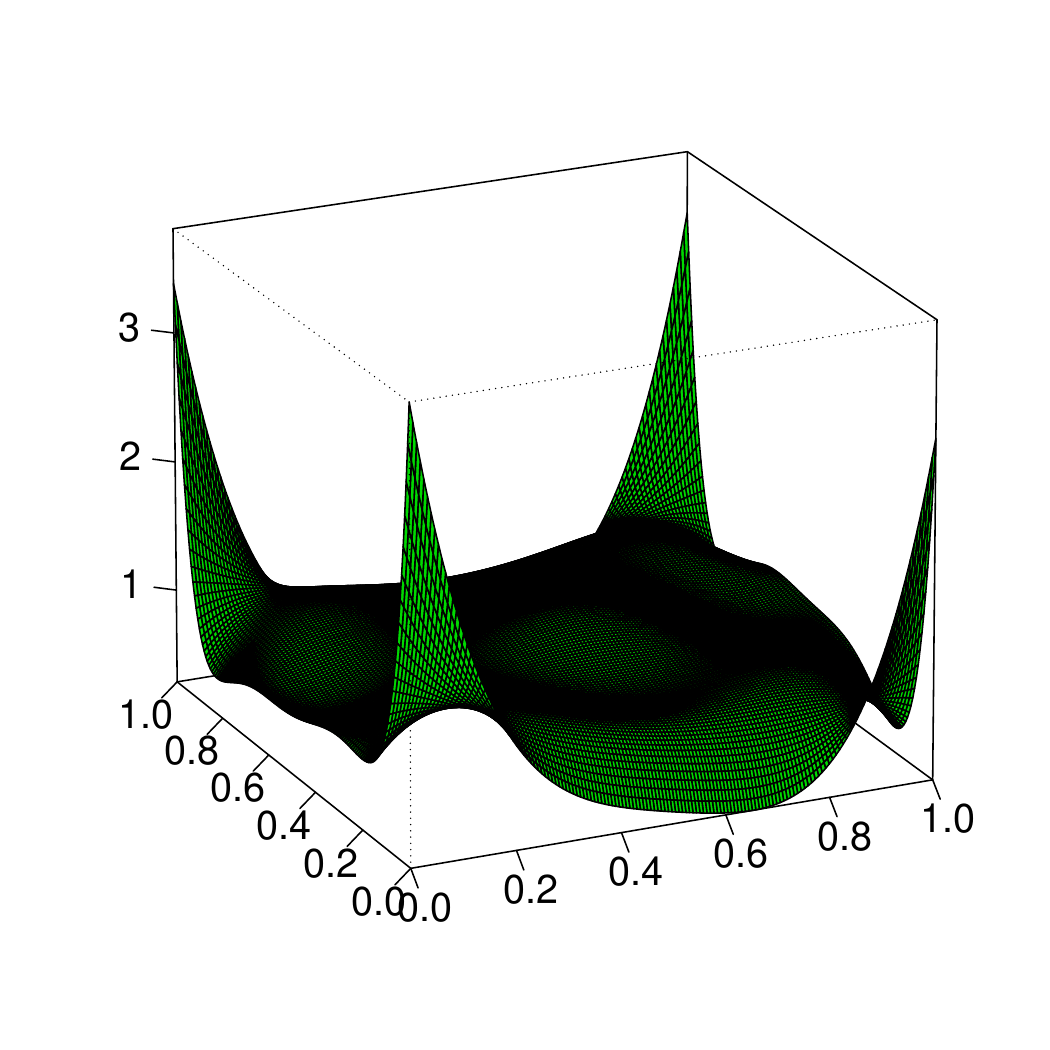}
  \includegraphics[height=.45\textheight,width=.55\textwidth,keepaspectratio,trim=.5cm .5cm .5cm .5cm]{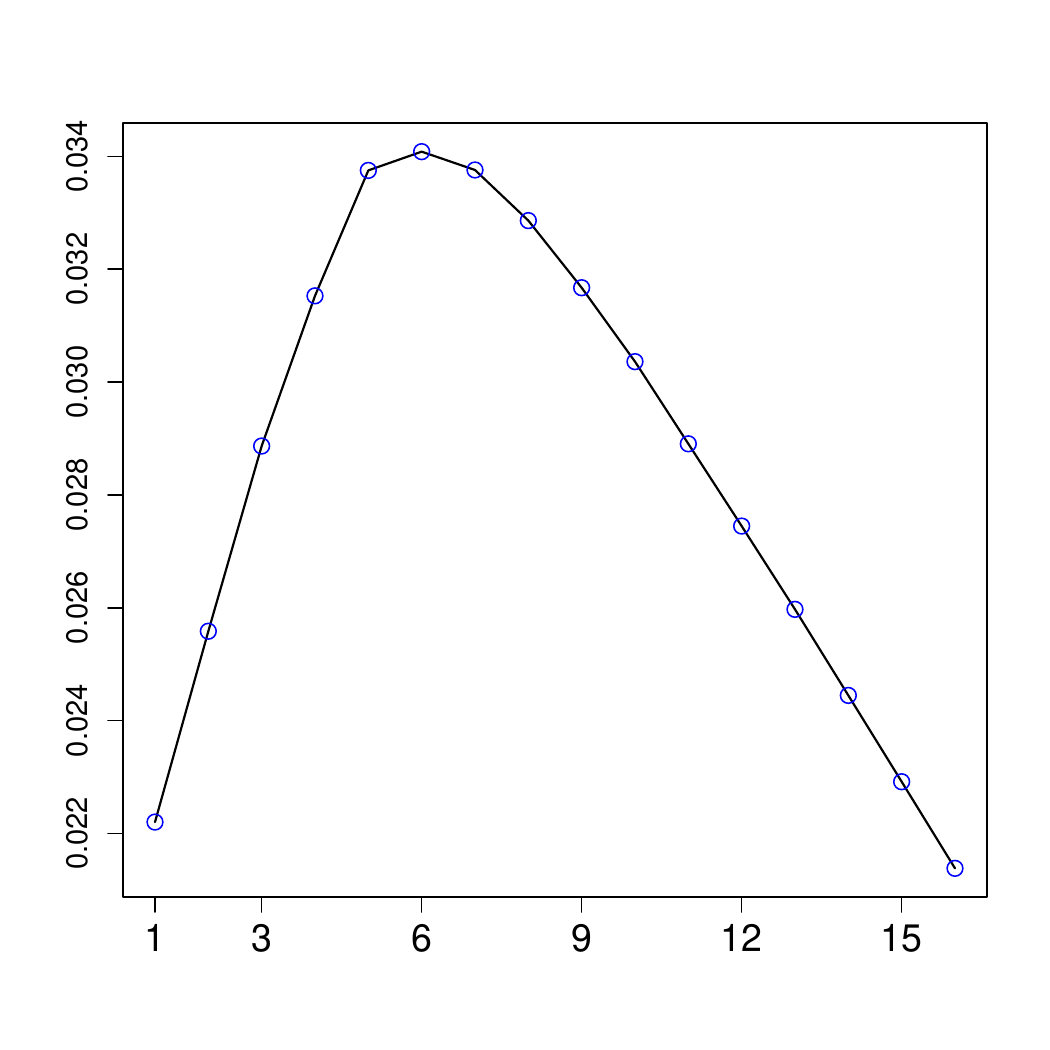}
\vspace{-1.4em}
\caption{Top: Nonparametric smooth serial copula density (lag $1$) estimate of S\&P return data. Bottom: BIC plot to select the significant LP-comoments computed in Eq. \ref{eq:lpc}.} \label{fig:LPtcop}
\end{figure*}

Now we apply this nonparametric copula estimation theory to model the temporal dependence structure of S\&P return data. The copula density estimate $\widehat{\cop}(u,v;Y(t),Y(t+1))$ based on the smooth LP comoments is displayed in Fig \ref{fig:LPtcop}. The shape of the copula density shows strong evidence of \emph{asymmetric tail-dependence}. Note that the \emph{dependence is only present in the extreme quantiles} - another well-known stylized fact of economic and financial time series.

\vskip.15em
\subsubsection{LP-Comomemnt of lag $h$}
Here we will introduce the concept of LP-comoment to get a complete understanding of the nature of serial dependence present in the data.
LP-comoment of lag $h$ is defined as the joint covariance of $\rm{Vec}(\YS)(t)$ and $\rm{Vec}(\YS)(t+h)$.

The lag 1 LP-comoment matrix for S\&P 500 return data is displayed below,
\beq \label{eq:lpc}
\LP\big[ Y(t),Y(t+1) \big]~=~\begin{bmatrix}
                              0.0705^{*} &-0.0617^{*} &0.0199  &0.0113\\
                              0.0074  &0.1542^{*} &0.0077  &0.0652^{*}\\
                             -0.0104 &-0.0071 &0.0262 &-0.0355\\
                              0.0166 &0.0438^{*} &0.0113  &0.0698^{*}
                             \end{bmatrix} \eeq
To identify the significant elements we first rank order the squared LP comoments. Then we take the penalized cumulative sum of $m$ comoments using BIC criterion $2m\log(n)/n$, $n$ is sample size and choose the $m$ for which BIC is maximum. The complete BIC path for S\&P 500 data is shown in figure \ref{fig:LPtcop}, which selects top six comments also denoted by ${*}$ in the LP-comoment matrix display (Eq. \ref{eq:lpc}). By making all those uninteresting ``small'' comoments equal to zero we get the ``smooth'' LP Comoment matrix denoted by $\hLP$ . The linear auto-correlation is capture by the $\LP[1,1;Y(t),Y(t+1)] = \Ex[\YS_1(t) \YS_1(t+1)]$ term. The presence of higher order significant terms in the LP comoment matrix indicate the possible \emph{nonlinearity}. Another interesting point to note that $\corr[Y(t),Y(t+1)]=.027$, where as the auto-correlation between the mid-rank transformed data $\corr[\Fm(Y(t)),\Fm(Y(t+1))] = .071$, considerably larger and picked by BIC criterion. This is an interesting fact as it indicates rank-transform time series ($\YS_1(t)$) is \emph{much more predictable} than the original raw time series $Y(t)$.

\subsubsection{LP-Correlogram, Evidence and Source of Nonlinearity}
We provide nonparametric exploratory test for  (non)linearity (spectral domain test is given in Sec \ref{sec:spec}). Plot the correlogram of $\YS_1(t), \ldots, \YS_4(t)$: (a) diagnose possible nonlinearity, and (b) identify possible sources. This constitutes the important building block for methods of model identification. LP-correlogram generalizes the classical sample autocorrelation function (ACF). Applying the \texttt{acf()} R function  on $\rm{Vec}(\YS)(t)$ generates the graphical display of our proposed LP-correlogram plot.
\begin{figure*}[!htb]
 \centering
 \includegraphics[height=.5\textheight,width=\textwidth,keepaspectratio,trim=.5cm .5cm .5cm .5cm]{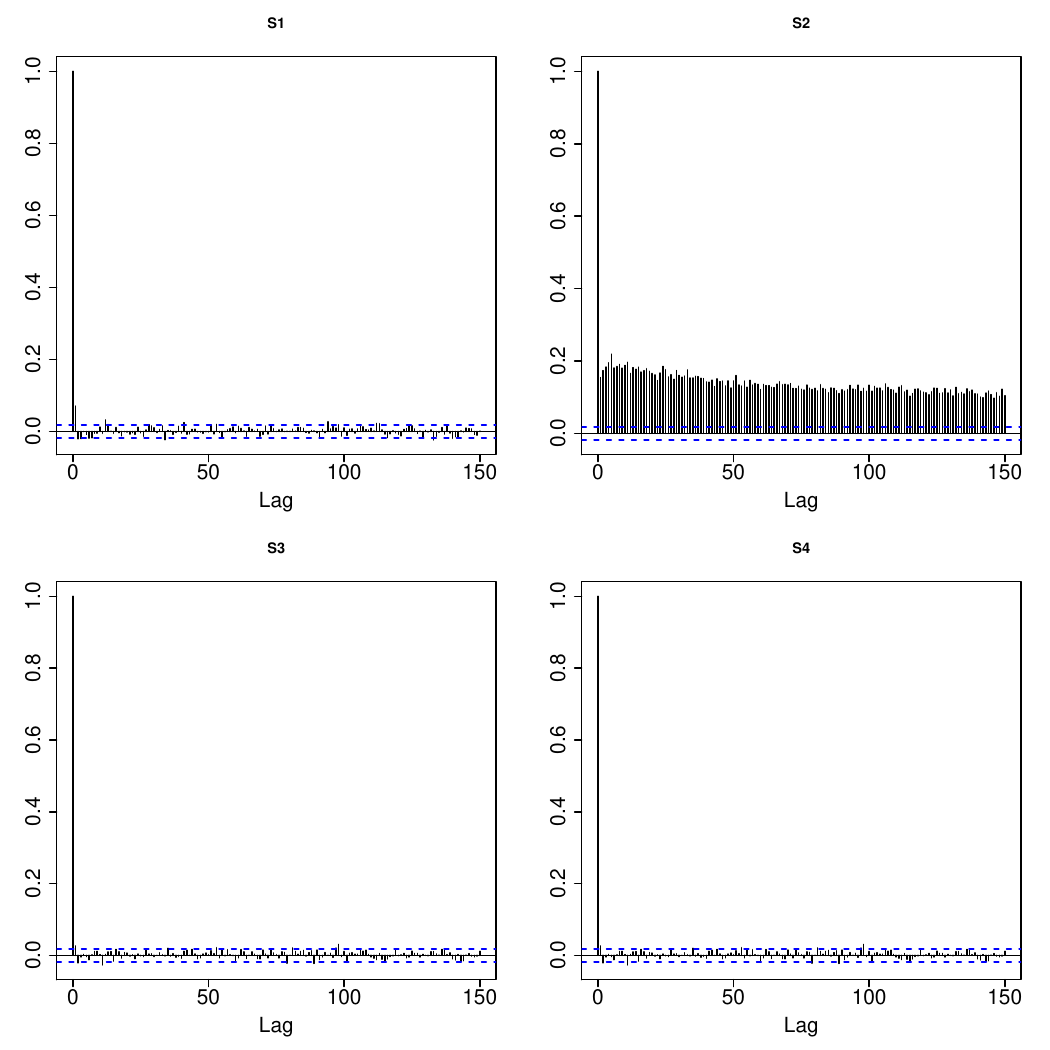} \\
\vspace{-.4em}
\caption{LP-correlogram: sample autocorrelations of LP-Transformed time series. The decay rate of the sample autocorrelations of $\YS_2(t)$ appears to be much slower than the exponential decay of ARMA process, implying possible long memory behaviour.} \label{fig:LPtacf}
\end{figure*}
\vskip.15em

Fig \ref{fig:LPtacf} shows the LP-correlogram of S\&P stock return data. Panel A shows the absence of linear autocorrelation which is known as a \emph{efficient market hypothesis} in finance literature. A prominent auto-correlation pattern for the series $\YS_2(t)$ (panel B Fig \ref{fig:LPtacf}) is the source of nonlinearity. This fact is known as \emph{``volatility clustering''}, which says that large price fluctuation is more likely to be followed by large price fluctuations. Also the slow decay of the autocorrelation of the series $\YS_2(t)$ can be interpreted as an indication of \emph{long-memory volatility structure}.

\subsubsection{AutoLPinfor: Nonlinear Correlation Measure}\label{sec:lpi}
We display the sample AutoLPinfor plot - diagnostic tool for nonlinear autocorrelation. We define the lag $h$ AutoLPinfor as the squared Frobenius norm of the smooth-LP comoment matrix of lag $h$,
\beq
\rm{AutoLPinfor}(h)~=~\sum_{j,k}\big|\,   \LP[j,k;Y(t),Y(t+h)]   \, \big|^2,
\eeq
where sum is over BIC selected $j,k$ for which LP comoments are significantly non-zero.

\begin{figure*}[t]
 \centering
 \includegraphics[height=.3\textheight,width=\textwidth,keepaspectratio,trim=.7cm .5cm .7cm .9cm]{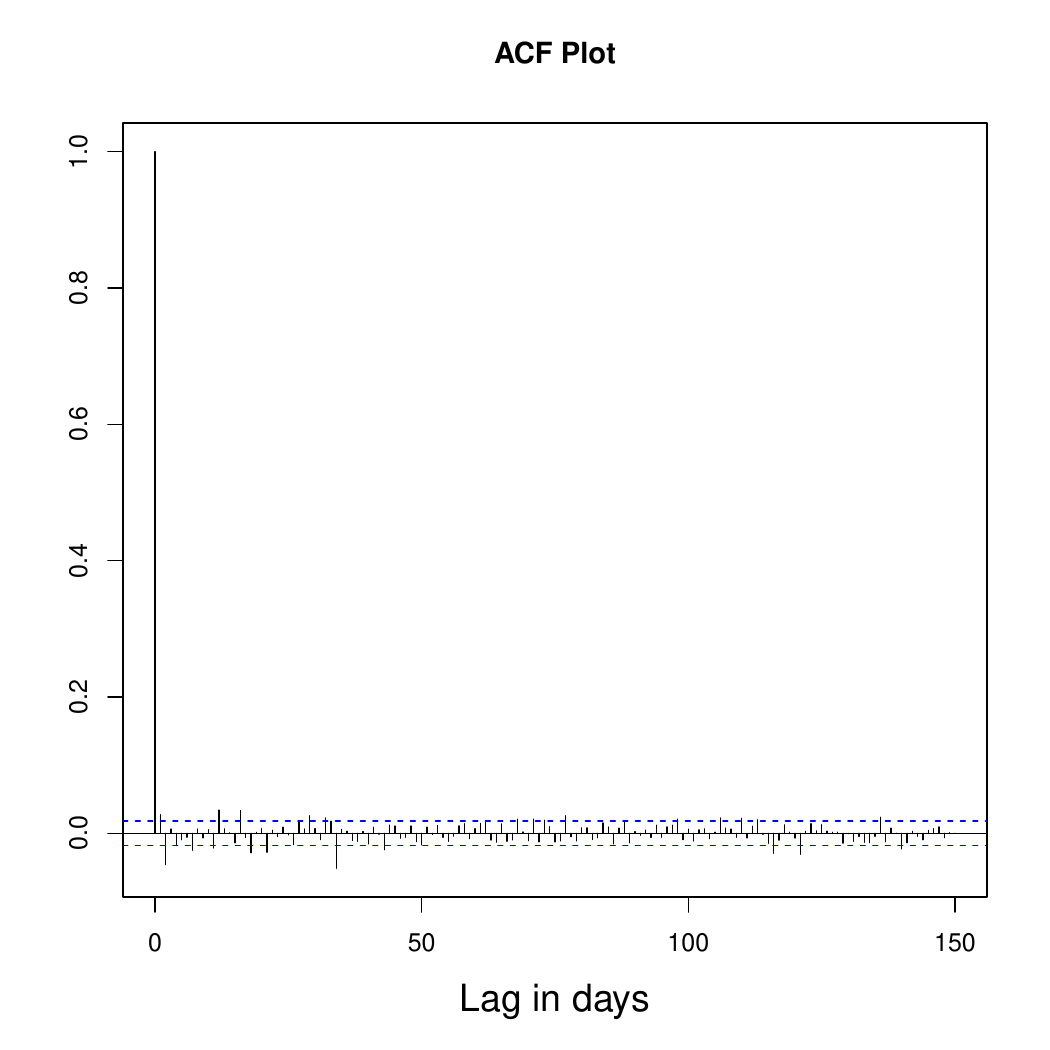}~~~~
 \includegraphics[height=.3\textheight,width=\textwidth,keepaspectratio,trim=.7cm .5cm .7cm .5cm]{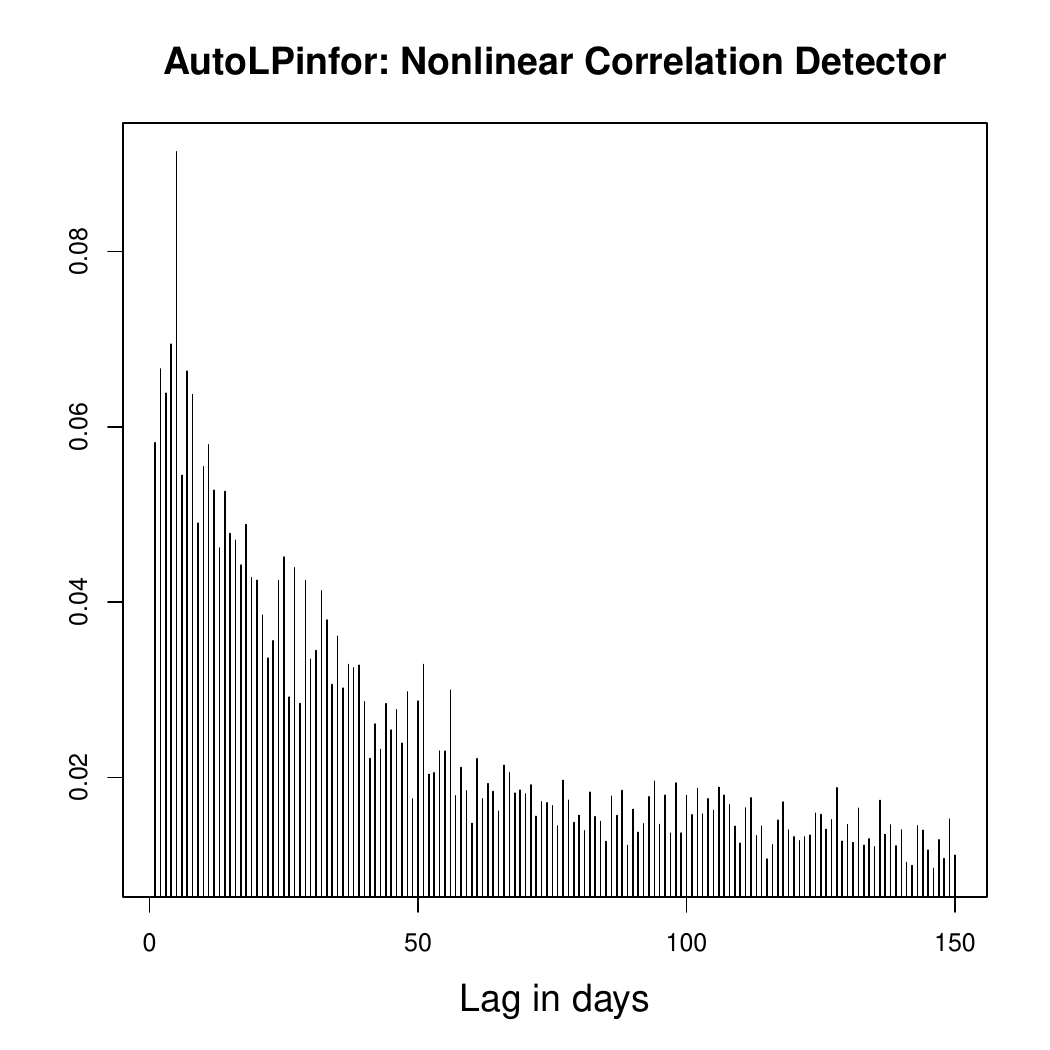}
\\
\vspace{-.4em}
\caption{Left: ACF plot of S\&P 500 data; Right: AutoLPinfor Plot upto lag $150$.} \label{fig:LPtauto}
\end{figure*}

Our robust nonparametric measure can be viewed as capturing the deviation of copula density from uniformity:
\bea
\rm{AutoLPinfor}(h)~=~ \iint \cop^2[u,v;Y(t),Y(t+h)]\dd u \dd v\,-\,1,
\eea
which is closely related to the entropy measure of association proposed in \cite{grangerlin}
\bea
\rm{Granger}\textrm{--}\rm{Lin}(h)=\iint \cop[u,v;Y(t),Y(t+h)]\, \log \cop[u,v;Y(t),Y(t+h)] \dd u \dd v.
\eea
It can be showed using Taylor series expansion that asymptotically
\beq
\rm{AutoLPinfor}(h)~\approx~ 2 \times \rm{Granger}\textrm{--}\rm{Lin}(h).
\eeq
An excellent discussion on the role of information theory methods for unified time series analysis is given in \cite{parzen92b} and \cite{brillinger04}. For an extensive survey of tests of independence for nonlinear processes see Chapter 7.7 of \cite{terasvirta2010book}. AutoLPinfor is a \emph{new} information theoretic nonlinear autocorrelation measure which detects \emph{generic} association and serial dependence present in a time-series. Contrast the AutoLPinfor plot for S\&P 500 return data shown in Fig \ref{fig:LPtauto} with the acf plot (left panel). This underlies the need for building a nonlinear time series model, which we will be discussing next.

\subsubsection{Nonparametric Estimation of Blomqvist's Beta}
Estimate the Blomqvist's $\be$ (also known as medial correlation coefficient) of lag $h$ by using the LP-copula estimate in the following equation,
\beq
\widehat{\be}_{\rm{LP}}(h;Y(t))~:=~-1 \,+\, 4 \int_{0}^{1/2}\int_{0}^{1/2}\, \widehat{\cop}\big[u,v;Y(t),Y(t+h)   \big] \dd u \dd v
\eeq

The $\be$ values $-1, 0$ and $1$ interpreted as reverse correlation, independence and perfect correlation, respectively. Note that,
\beas
\mbox{Blomqvist's}\,\be&&: \mbox{Normalized distance of copula distribution} \Cop(u,v)\, \mbox{from independence copula} \,uv \\
\mbox{AutoLPinfor}&&: \mbox{Distance of copula density} \cop(u,v)\, \mbox{from uniformity} \,1 .
\eeas
For S\&P 500 return data we compute the following dependence numbers,
\begin{figure*}[t]
 \centering
 \includegraphics[height=.35\textheight,width=\textwidth,keepaspectratio,trim=2.6cm 1.2cm 2.6cm 1.4cm]{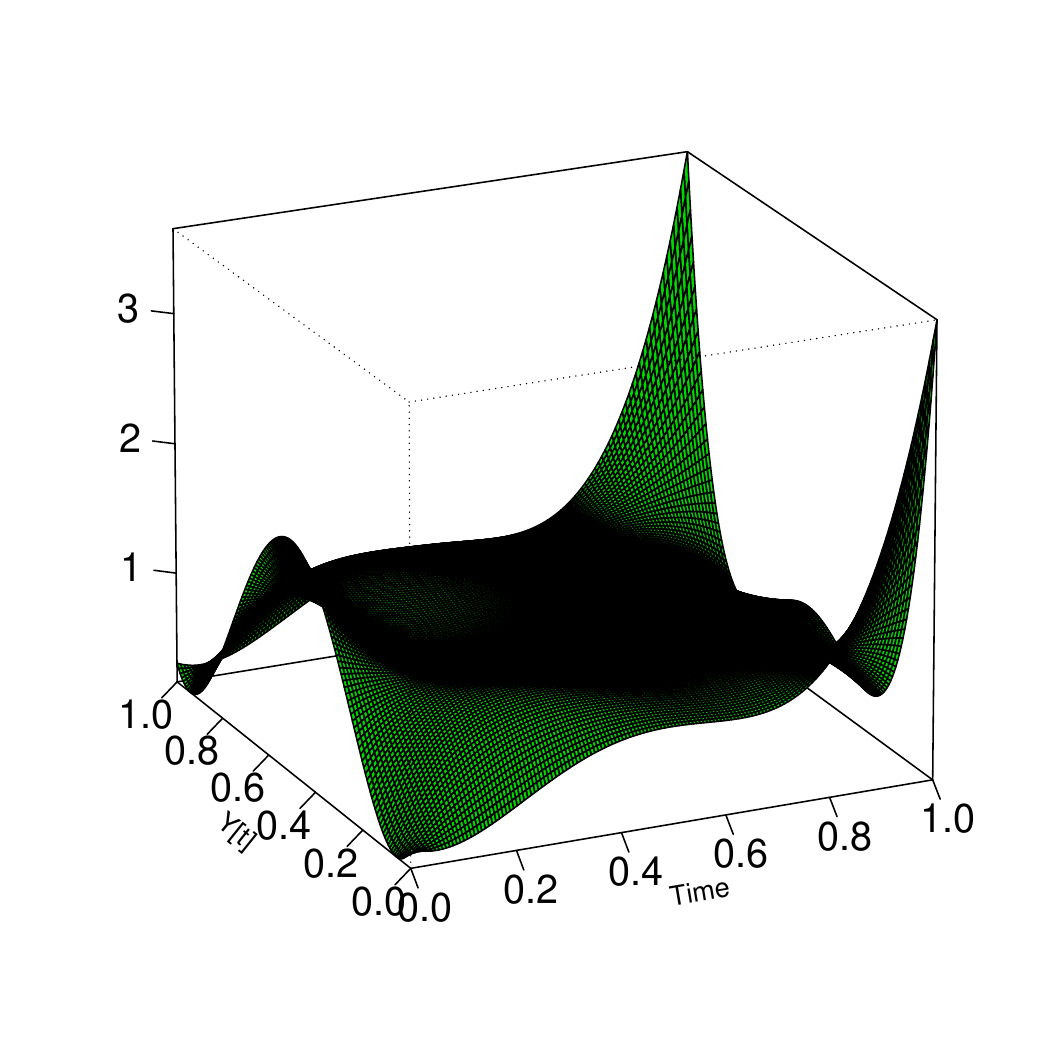} \\
\vspace{-1.4em}
\caption{LP copula diagnostic for detecting non-stationarity in S\&P 500 return data.} \label{fig:LPtnon}
\end{figure*}
\beas
&&\widehat{\be}_{\rm{LP}}(1;Y(t))=\,0.0528 \\
&&\widehat{\be}_{\rm{LP}}(1;\YS_1(t))=\,0.0528\\
&&\widehat{\be}_{\rm{LP}}(1;\YS_2(t))=\,0.0729\\
&&\widehat{\be}_{\rm{LP}}(1;\YS_3(t))=\,0.0\\
&&\widehat{\be}_{\rm{LP}}(1;\YS_4(t))=\,0.003.
\vspace{-.7em}
\eeas
\subsubsection{Nonstationarity Diagnosis, LP-Comoment Approach}
Viewing the time index $T=1,\ldots,n$ as covariate we propose a nonstationarity diagnosis based on LP-comoments of $Y(t)$ and the time index variable $T$. Our treatment has the ability to detect time varying nature of mean, variance, skewness and so on represented by various custom-made LP-transformed time series.

For S\&P data we computed the following LP comoment matrix to investigate the nonstationarity:
\beq \label{eq:lpnon}
\LP\big[ T,Y(t) \big]~=~\begin{bmatrix}
                              0.012  &0.180^{*} &-0.010  &0.058^{*}\\
                              -0.005 &-0.034 &-0.036  &0.080^{*}\\
                             -0.016  &0.115^{*}  &0.001 &-0.001\\
                              0.024 &-0.040 &-0.010  &0.049^{*}
                             \end{bmatrix} \eeq

This indicates presence of \emph{slight non-stationarity behavior of variance or volatility} ($YS_2(t)$) and \emph{kurtosis} of tail-thickness($YS_4(t)$). Similar to AutoLPinfor we propose the following statistic for detecting nonstationarity
\beq
\rm{LPinfor}[Y(t),T]=~\sum_{j,k}\big|\,   \LP[j,k;T,Y(t)]   \, \big|^2.
\eeq
We can also generate the corresponding smooth copula density of $(T,Y(t))$ based on smooth $\LP\big[ T,Y(t) \big]$ matrix to visualize the time varying information as in Fig \ref{fig:LPtnon}.

\subsection{Local Dependence Modeling}
\subsubsection{Quantile Correlation Plot and Test for Asymmetry}
Display the quantile correlation plot, a copula distribution based graphical diagnostic to \emph{visually examine the asymmetry of dependence}. The goal is to get more insight into the \emph{nature of tail-correlation}.

\begin{figure*}[!htb]
 \centering
 \includegraphics[height=.4\textheight,width=\textwidth,keepaspectratio,trim=.8cm 1.2cm .8cm .7cm]{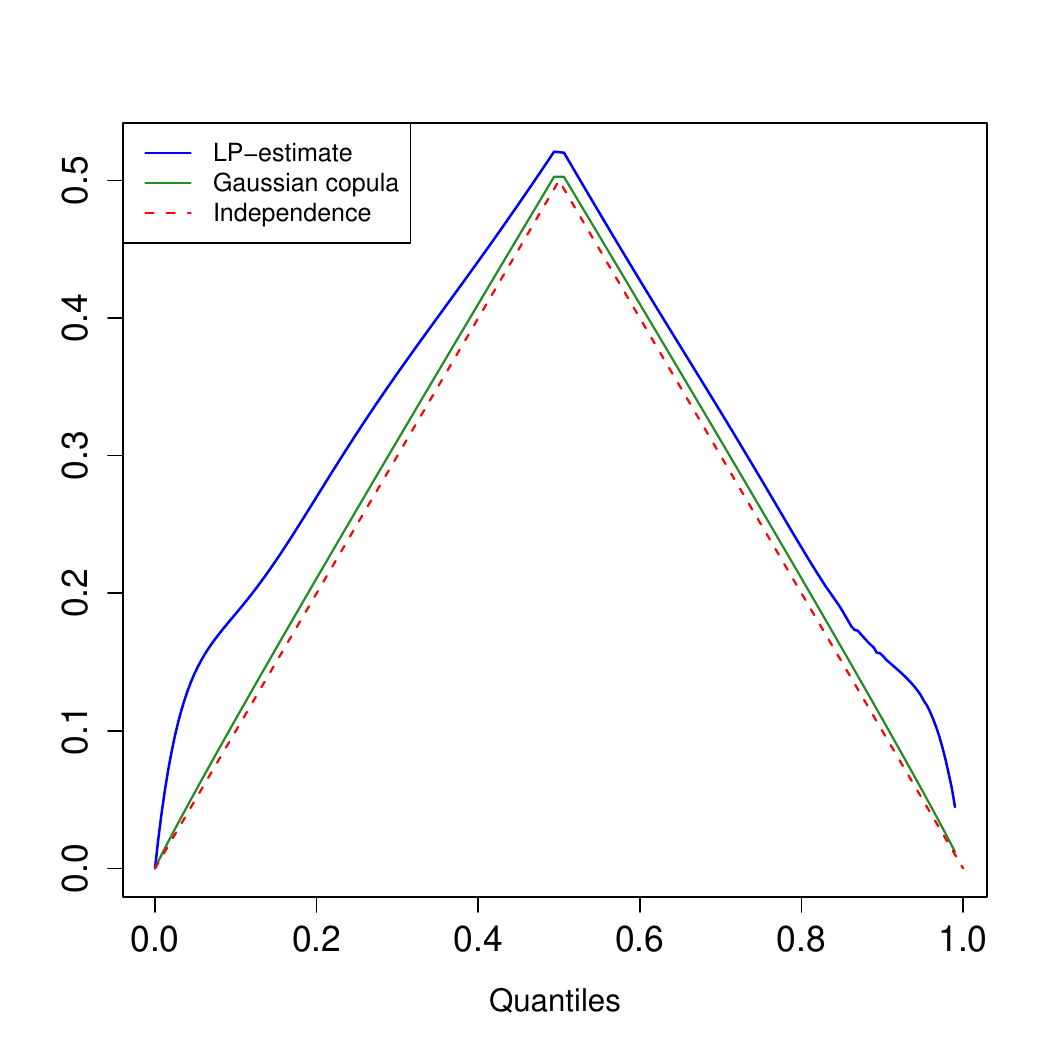} \\
\caption{Estimated Quantile Correlation Function (\texttt{QCF}) $\widehat{\la}_{\rm{LP}}[u;Y(t),Y(t+1)]$. Detects asymmetry in the tail-dependence between lower-left quadrant and upper-right quadrant for S\&P 500 return data. The red dotted line denotes the quantile correlation function under dependence. The dark-green line shows the quantile correlation curve for fitted Gaussian copula.} \label{fig:LPttail}
\end{figure*}

Motivated by the concept of lower and upper tail-dependence coefficient we define the \emph{Quantile Correlation Function} (\texttt{QCF}) as the following in terms of copula distribution function of $(Y(t),Y(t+h))$ denoted by $\Cop(u,v;Y(t),Y(t+h)):=\Cop(u,v;h)$,
\beq
\la\big[u;Y(t),Y(t+h)\big]\,:=\, \dfrac{\Cop(u,u;h)}{u}\, \ind_{\big\{u \le .5 \big\}} ~+~ \dfrac{1 -2u +\Cop(u,u;h)}{1-u} \, \ind_{\big\{ u > .5\big\} },
\eeq
Our nonparametric estimate of quantile correlation function is based on the LP-copula density which we denote as $\widehat{\la}_{\rm{LP}}\big[u;Y(t),Y(t+h)\big]$. Fig \ref{fig:LPttail} shows the corresponding quantile correlation plot for S\&P 500 data. The dotted-line represent \texttt{QCF} under Independence assumption. Deviation from this line helps us to better understand the nature of asymmetry. We compute $\widehat{\la}_{\rm{G}}[u;Y(t),Y(t+h)]$ using the fitted Gaussian copula
\beq \widehat{\Cop}_{\rm{G}}(u,v;Y(t),Y(t+h))~=~\Phi(\Phi^{-1}(u),\Phi^{-1}(v);\hat{\Si}=S)\eeq
where $S$ is sample covariance matrix. The dark-green line in Fig \ref{fig:LPttail} shows the corresponding curve which is almost identical with the ``no dependence'' curve, albeit misleading. The reason is Gaussian copula is characterized by linear correlation, while S\&P data is highly nonlinear in nature. As the linear auto-correlation of stock return is almost zero, we have approximately $\Phi(\Phi^{-1}(u),\Phi^{-1}(u);\hat{\Si}=S) \approx \Phi(\Phi^{-1}(u))\Phi(\Phi^{-1}(u)) = u^2$. Similar to Gaussian copula there are several others parametric copula families, which can give similar misleading conclusions. This simple illustration reminds us the pernicious effect of \textit{not} ``looking into the data''.
\vskip.15em

\subsubsection{Conditional LPinfor Dependence Measure}
For more transparent and clear insights into the asymmetric nature of the tail dependence we need to introduce the concept of \emph{conditional dependence}. In what follows we propose a conditional LPinfor function $\LPinfor(Y(t+h)|Y(t)=Q(u;Y(t)))$ - quantile based diagnostic for tracking how the dependence of $Y(t+h)$ on $Y(t)$ changing at various quantiles.

To quantify the conditional dependence, we seek to estimate $f(y;Y(t+h)|Y(t))/f(y;Y(t+h))$. A brute force approach estimates separately the conditional distribution and the unconditional distribution and take the ratio to estimate this arbitrary function. An alternative elegant way is to recognize that by ``going to quantile domain'' (i.e., $Y(t+h)=Q(v;Y(t+v))$ and $Y(t)=Q(u;Y(t))$) we can interpret the ratio as ``slices'' of copula density, which we call conditional comparison density:
\beq \label{eq:ccd}
d\big[v;Y(t+h),Y(t+h)|Y(t)=Q(u;Y(t))\big]= 1 + \sum_j \LP[j;h,u]\, S_j(v;Y(t+h)),
\eeq
where the LP-Fourier orthogonal coefficients $\LP[j;h]$ are given by
\[\LP[j;h,u]~=~\sum_k \LP[j,k;Y(t),Y(t+h)] \,S_k(u;Y(t)).\]
Define conditional LPinfor as
\beq
\LPinfor\big[Y(t+h)|Y(t)=Q(u;Y(t))\big]~=~\sum_j \big|\LP[j;h,u]\big|^2.
\eeq



\begin{figure*}[!htb]
 \centering
 \includegraphics[height=\textheight,width=.32\textwidth,keepaspectratio,trim=.5cm .5cm .5cm 1.2cm]{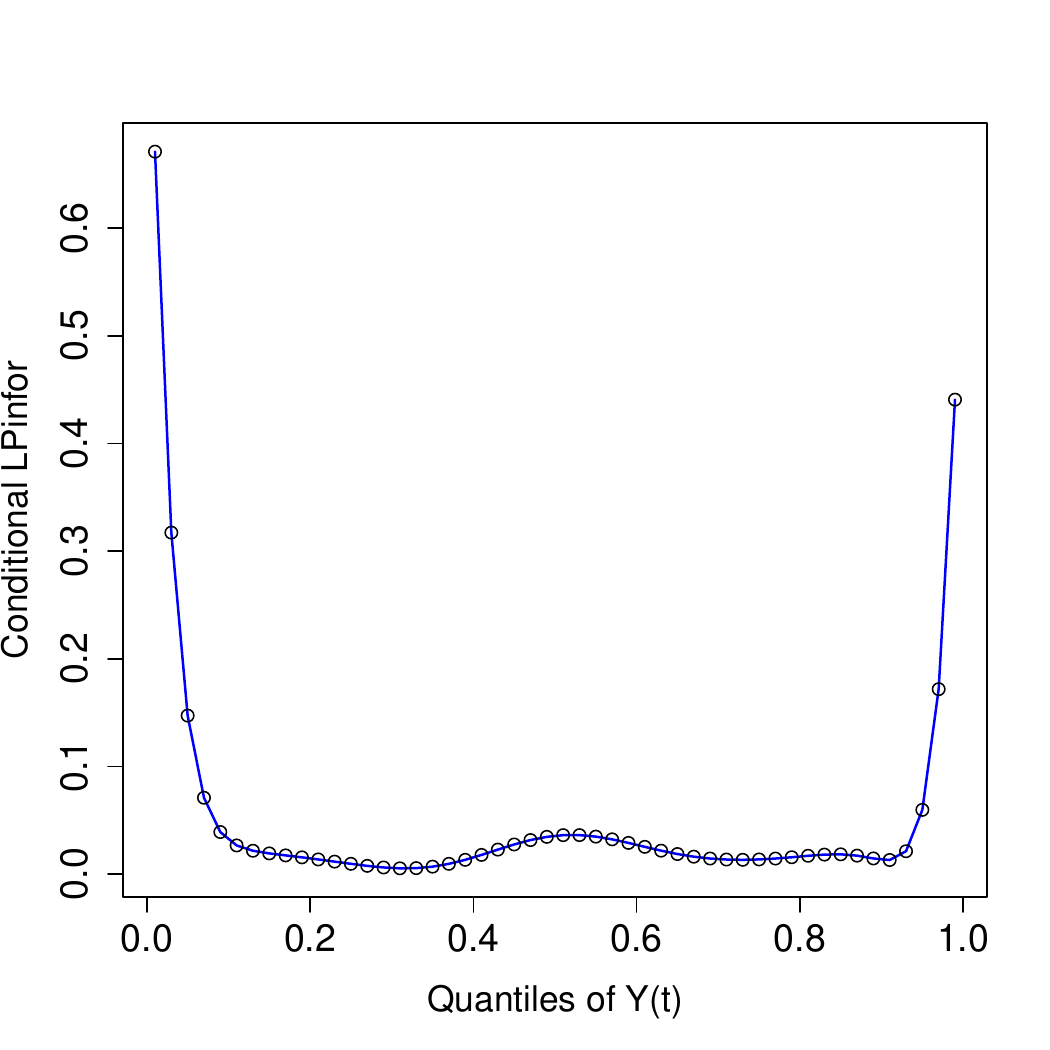}
  \includegraphics[height=\textheight,width=.32\textwidth,keepaspectratio,trim=.5cm .5cm .5cm .5cm]{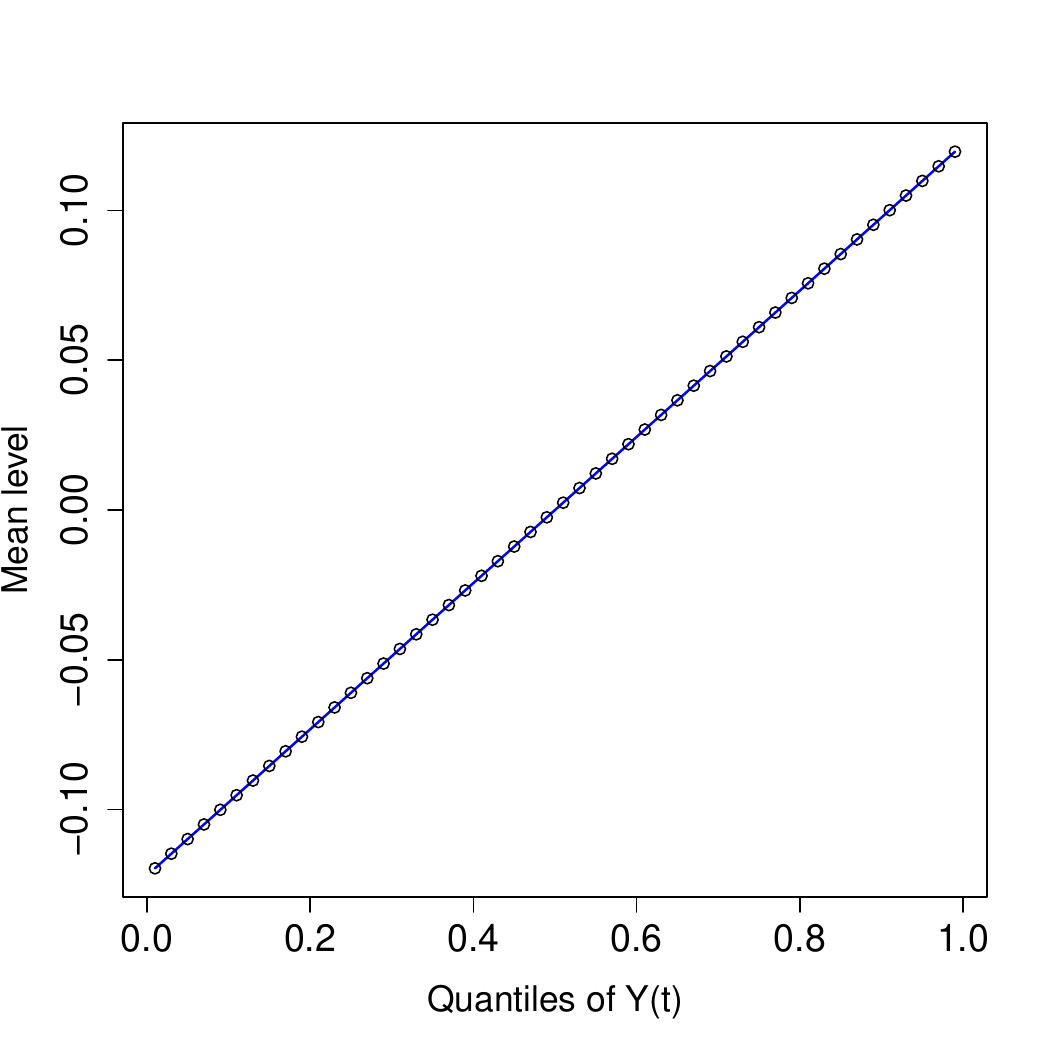}
  \includegraphics[height=\textheight,width=.32\textwidth,keepaspectratio,trim=.5cm .5cm .5cm .5cm]{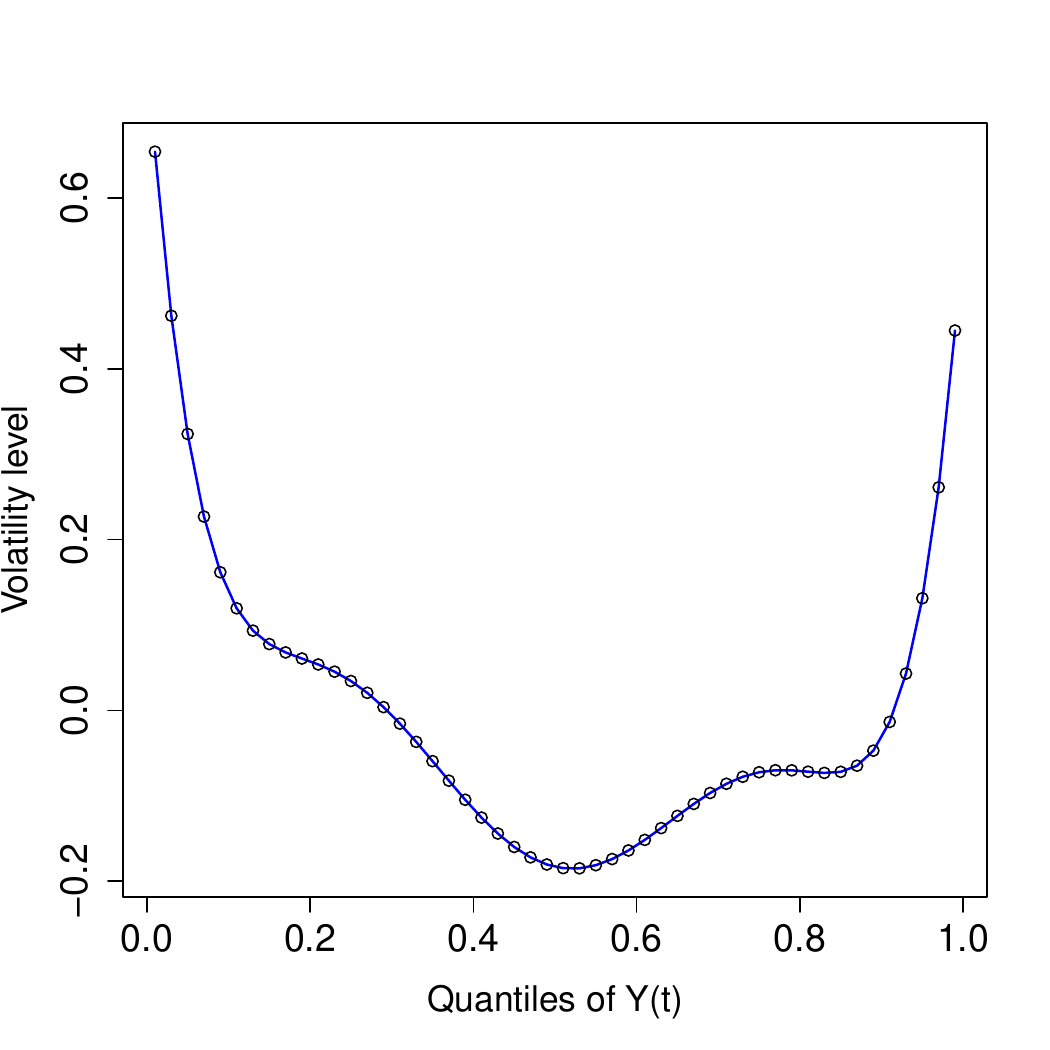}
\vspace{-.5em}
\caption{(a) The conditional LPinfor curve is shown for the pair $[Y(t),Y(t+1)]$. Clearly shows the asymmetric dependence in the tails and almost nothing is going on in between. (b) and (c) display how the mean and volatility levels of conditional distribution $f[y;Y(t+1)|Y(t)=Q(u;Y(t))]$ changes with respect to the unconditional marginal distribution $f(y;Y(t))$ at different quantiles.} \label{fig:cLPi}
\end{figure*}
\vskip.15em
%

We use this theory to investigate the conditional dependency structure of S\&P 500 return data. Fig \ref{fig:cLPi}(a) traces out the complete path of the estimated $\LPinfor[Y(t+h)\mid Y(t)=Q(u;Y(t)]$ function, which indicates the \emph{high asymmetric tail-correlation}. This conditional correlation curves can be viewed as ``local'' dependence measure.  An excellent discussion on this topic is given in Section 3.3.8 of \cite*{terasvirta2010book}.

At this point we can legitimately ask \emph{what aspects of the conditional distributions are changing most}? Fig \ref{fig:cLPi} (b,c) display only the two coefficients $\LP[1;h,u]$ and $\LP[2;h,u]$ for the S\&P 500 return data for the pairs $(Y(t),Y(t+1))$. These two coefficients represent how the mean and the volatility levels of the conditional density changing with respect to the unconditional reference distribution. The typical \emph{asymmetric shape of conditional volatility} shown in right panel of Fig \ref{fig:cLPi} (b,c) indicates what is known as the \emph{``leverage effect''} - future stock volatility negatively correlated with past stock return, i.e. stock volatility tends to increase when stock prices drop.

\subsection{Non-Crossing Conditional Quantile Modeling}
Display the nonparametrically estimated conditional quantile curves of $Y(t+h)$ given $Y(t)$. Our new modeling approach uses the estimated  conditional comparison density $\dhat(v;h,u)$ to simulate from $F[y;Y(t+h)|Y(t)=Q(u;Y(t))]$ by utilizing the given sample $\widetilde Q(u;Y(t))$ via accept-reject rule to arrive at the ``smooth'' nonparametric model for $\widehat{Q}[v;Y(t+h)|Y(t)=Q(u;Y(t)]$. See \cite{D13USA} for details about the method. Our proposed algorithm generates ``large'' additional simulated samples from the conditional distribution, which allows us to accurately estimate the conditional quantiles (especially the extreme quantiles). By construction, our method guaranteed to produce non-crossing quantile curves - thus tackles a challenging practical problem.

\begin{figure*}[!htb]
 \centering
 \includegraphics[height=.55\textheight,width=\textwidth,keepaspectratio,trim=1.2cm .7cm 1.2cm .5cm]{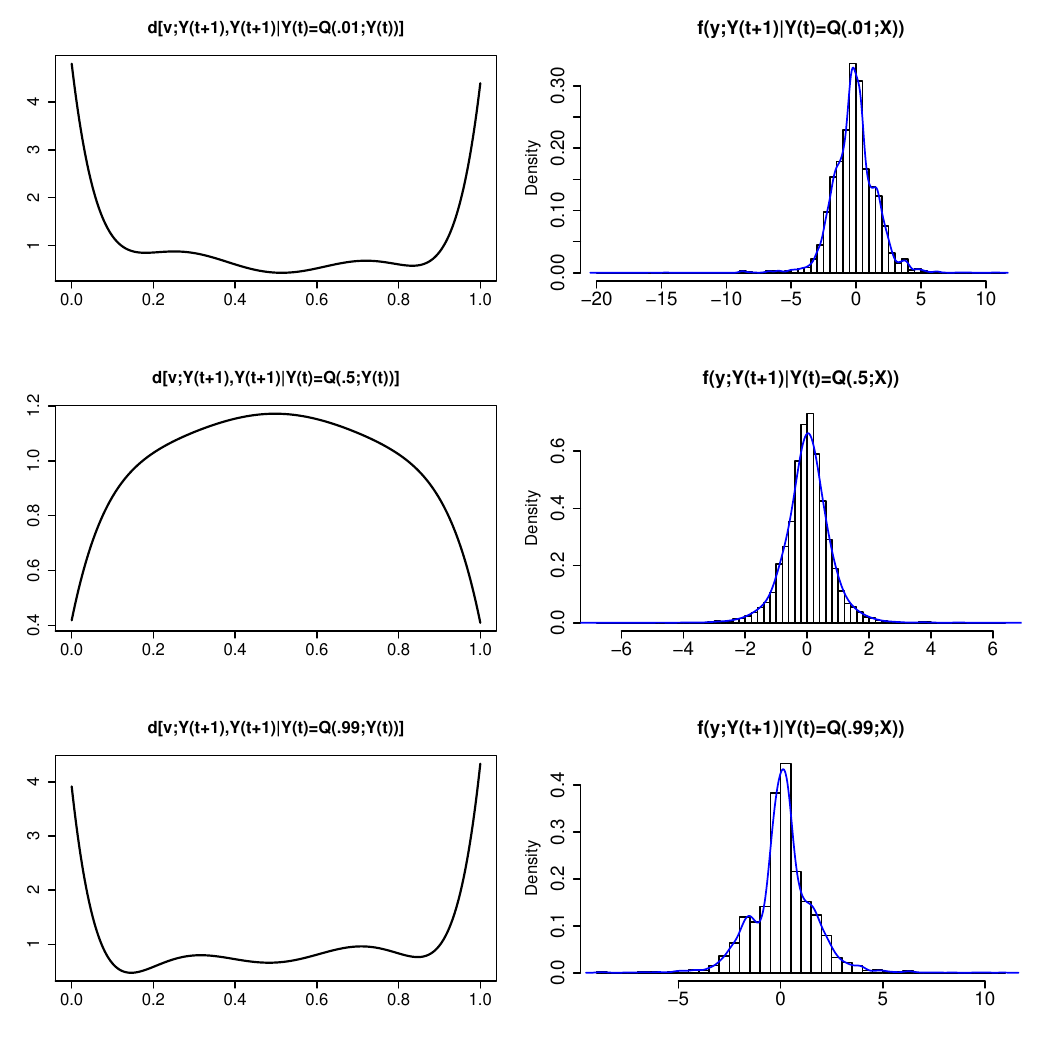} \\
\caption{Each row display the estimated conditional comparison density and the corresponding conditional distribution for $u=.01,.5,.99$.} \label{fig:cdist}
\end{figure*}
\vskip.15em

For S\&P 500 data we first nonparametrically estimate the conditional comparison densities $\dhat(v;h,u)$ shown in Left panel of Fig \ref{fig:cdist} for $F(y;Y(t))=.01,.5$ and $.99$, which can be thought of as a ``weighting function'' for unconditional marginal distribution to produce the conditional distributions:
\beq
\fhat\big[y;Y(t+h)|Y(t)=Q(u;Y(t))\big]~=~f(y;Y(t)) \, \times \, \dhat\big[ F(y;Y(t+h));h,u \big].
\eeq
This density estimation technique belongs to the skew-G modeling class \citep{deep16LSSD}. We simulate $n=10,000$ samples from the $\fhat(y;Y(t+h)|Y(t))$ by accept-reject sampling from $\dhat(v;h,u)$, $u=\{.01,.5,.99\}$. The histograms and the smooth conditional densities are shown in the right panel of Fig \ref{fig:cdist}. It shows some typical shapes in terms of long-tailedness.
\begin{figure*}[!htb]
 \centering
 \includegraphics[height=.35\textheight,width=\textwidth,keepaspectratio,trim=1.5cm .4cm 1.5cm 2cm]{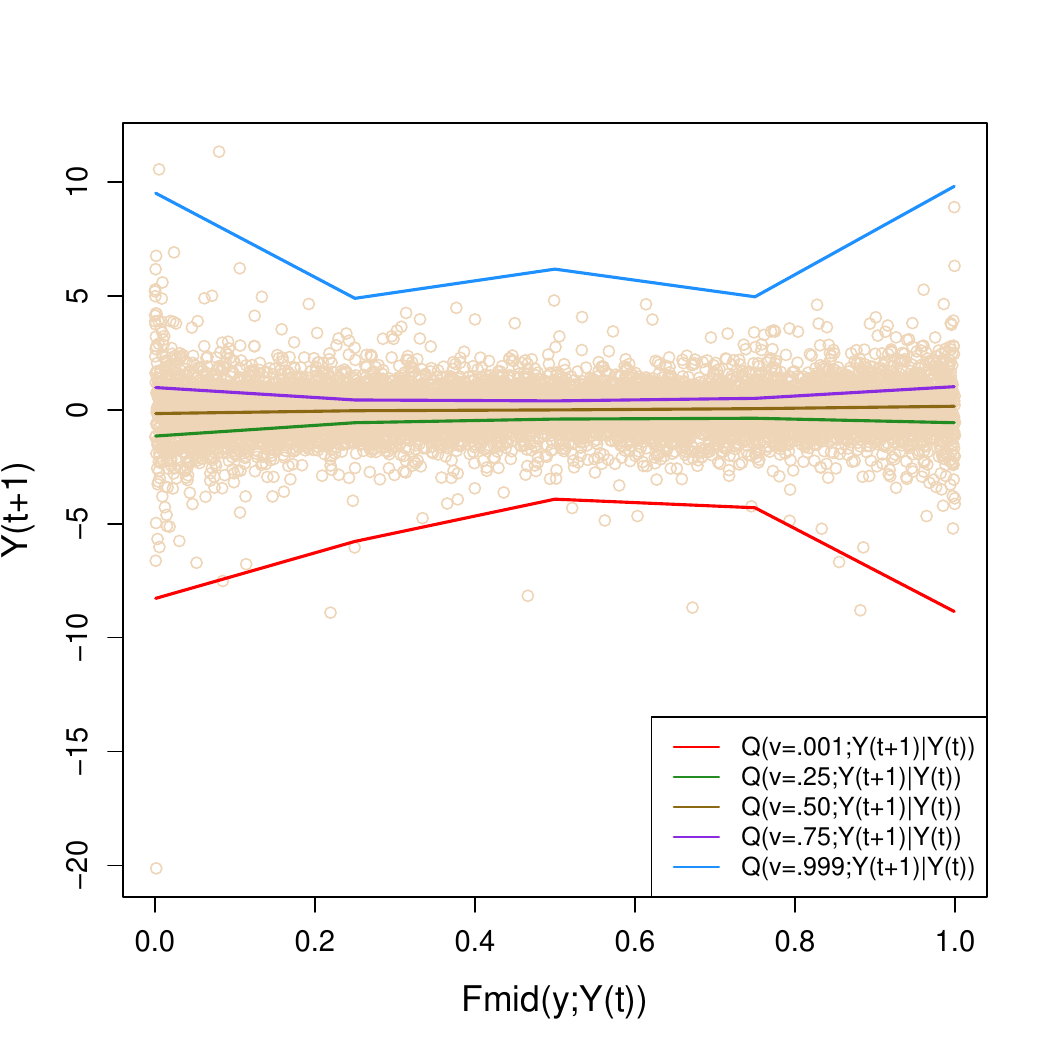}\\
\caption{The figure shows estimated non-parametric conditional quantile curves for S\&p 500 return data. The red solid line which represent $\widehat{Q}(.001;Y(t+1)|Y(t))$ is popularly known as one day .1\% conditional value at risk measure (CoVaR). } \label{fig:cq}
\end{figure*}

Next we proceed to estimate the nonparametric conditional quantiles $\Qhat(v;Y(t+h)|Y(t))$, for $v=.001,.25,.5,.75,.999$ from the simulated data. Fig \ref{fig:cq} shows the estimated conditional quantiles. The extreme conditional quantiles have a special significance in the context of financial time series. They sometimes popularly known as \emph{Conditional Value at Risk (CoVaR) - currently the most popular quantitative risk management tool }(see \cite{covar,engle04}). The red solid line in Fig \ref{fig:cq} is the $\Qhat\big[.001;Y(t+1)\mid Y(t)=Q(u;Y(t))\big]$, which is known as .1\% CoVar function for one-day holding period for S\&P 500 daily return data. Although the upper conditional quantile curve $\Qhat(.999;Y(t+1)|Y(t))$ (blue solid line) show symmetric behaviour around $F(y;Y(t))=.5$, the lower quantile has a prominent asymmetric shape. These conditional quantiles give ultimate description of the auto-regressive dependence of S\&P 500 return movement in the tail-region.

\subsection{Nonlinear Spectrum Analysis} \label{sec:spec}
Here we extends the concept of spectral density for nonlinear processes.
Display the LPSpectrum - autoregressive (AR) spectral density estimates of $\YS_1(t),\ldots,\YS_4(t)$.
Spectral density for each LP-transformed series is defined as
\bea
f(\om;\YS_j) & = & \sum_{h} \LP[j,j;Y(t),Y(t+h)]\,e^{-i2\pi h \om},\qquad \,-1/2 < \om < 1/2 \nonumber \\
&=&\sum_{h}  \mbox{Cov}[\YS_j(t),\YS_j(t+h)] \,e^{-i2\pi h \om},\qquad -1/2 < \om < 1/2.
\eea

We separately fit the univariate AR model for the components of $\rm{Vec}(\YS)(t)$ and use BIC order selection criterion to select the ``best'' parsimonious parametrization using Burg method.

Finally we use the estimated model coefficients to produce the ``smooth'' estimate of spectral density function (see Eq. \ref{eq:arsdf}). The copula spectral density is defined as
\beq
f(\om;u,v)\,=\,\sum_{h} \cop(u,v;h)\,e^{-i2\pi h \om},\qquad \,-1/2 < \om < 1/2.
\eeq
To estimate the copula spectral density we use the LP comoment based nonparametric copula density estimate. Note that both the serial copula (3.12) and the corresponding spectral density (3.25) captures the same amount of information for serial dependence of $\{Y(t)\}$. For that reason, we recommend to compute AutoLPinfor as a general dependence measure for non-gaussian nonlinear processes. 
\begin{figure*}[!htb]
 \centering
 \includegraphics[height=.4\textheight,width=\textwidth,keepaspectratio,trim=1cm .5cm .5cm 1cm]{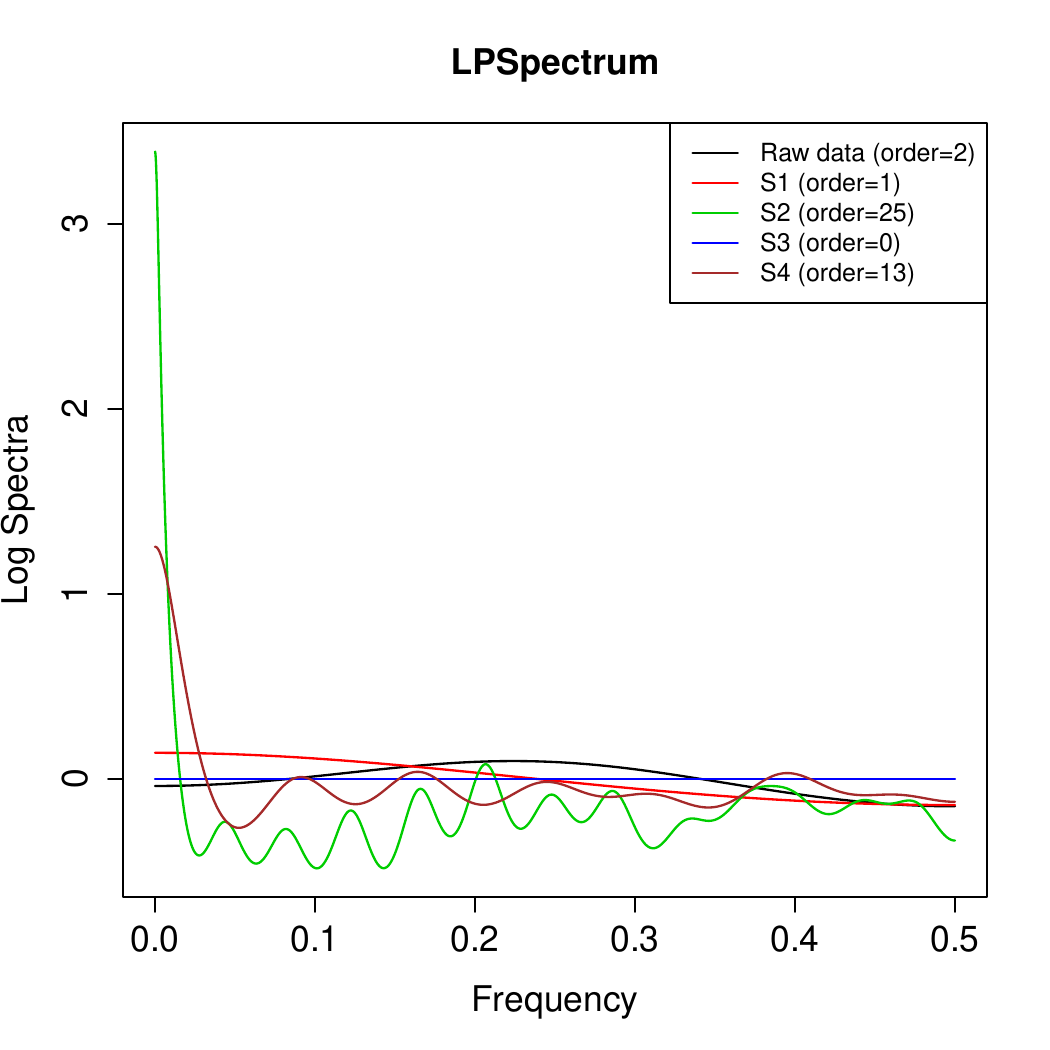} \\
\vspace{-.4em}
\caption{LPSpec: AR spectral density estimate for S\&P 500 return data. Order selected by BIC method. This provides a diagnostic tool for providing evidence of hidden periodicities in non-Gaussian nonlinear time series.} \label{fig:LPtspec}
\end{figure*}
Application of our LPSpectral tool on S\&P 500 return data is shown in Fig \ref{fig:LPtspec}. Few interesting observations: (i) the conventional spectral density (black solid line) provides no insight into the (complex) serial dependency present in the data. (ii) The nonlinearity in the series is captured by the interesting shapes of our specially-designed times series $\YS_2(t)$ and $\YS_4(t)$, which classical (linear) correlogram based spectra can not account for. (iii) The shape of the spectra of $\cZ(Y(t))$ and the rank-transformed time series $\YS_1(t)$ looks very similar. and a (iv) pronounced singularity near zero of the spectrum of $\YS_2(t)$ hints some kind of ``long memory'' behavior. This phenomena also known as regular
variation representation at frequency $\om=0$ \citep{granger1980LM}.

A quick diagnostic measure for screening significant spectrums can be computed via the information number $2 \int_0^{1/2} \log \widehat{f}(\om; S_j)\dd \om$. The LPSpectrum methodology is highly robust and thus can tackle the heavy-tailed S\&P data quite successfully.
\subsection{Nonparametric Model Specification} 

The ultimate goal of  empirical time series analysis is nonparametric model identification. To model the univariate stationary nonlinear process we specify multiple autoregressive model based on $\rm{Vec}(\YS)(t)\, =\,\big[ \YS_1(t), \ldots,   \YS_k(t) \big]^{T}$ of the form
\beq
\rm{Vec}(\YS)(t)~=~\sum_{k=1}^{m}  A(k;m) \,\rm{Vec}(\YS)(t-k)  + \bm \ep(t).
\eeq
where $\bm \ep(t)$ is multivariate mean zero Gaussian white noise with covariance $\Si_m$.
These system of equation \textit{jointly} describe the dynamics of the nonlinear process and how it evolves over time. We use BIC criterion to select the model order which selects the model order $m$ which minimizes
\beq
\rm{BIC(m)~=~\log|\widehat{\Si}_{m}|}\,+\, m k^2 \dfrac{\log T}{T}.
\eeq

We carry out this steps for our S\&P 500 return data. We estimate our multiple AR model based on $\rm{Vec}(\YS)(t)\, =\,\big[ \YS_1(t), \YS_2(t),\YS_4(t)  \big]^{T}$. We discard $\YS_3(t)$ due to it's flat spectrum (see Fig \ref{fig:LPtspec}). BIC selects ``best'' order $8$. Although, the complete description of the estimated model is clearly cumbersome, we provide below the approximate structure by selecting few large coefficients from the actual matrix equation. The goal is to interpret the coefficients (statistical parameters) of the estimated model and relate it with economic theory (scientific parameters/theory). This multiple AR LP-model (in short we call \texttt{LPVAR}) is given by
{\small
\bea
\YS_1(t)&\approx&  .071 \YS_1(t-1)- .024 \YS_1(t-2) + \ep_1(t)\nonumber \\
\YS_2(t)&\approx& -.063 \YS_1(t-1) -.075 \YS_1(t-2) + .06 \YS_2(t-2) + .123\YS_2(t-5)+ .04\YS_4(t-2) + \ep_2(t)\nonumber \\
\YS_4(t)&\approx&.04 \YS_4(t-1) + .038 \YS_4(t-2) + .04 \YS_2(t-3)+ \ep_4(t).
\eea}
and the residual covariance matrix is
\[\widehat{\Si_{8}} =\begin{bmatrix}
                              0.993 &-0.001 &-0.002\\
                              -0.001 &0.853 &-0.058\\
                              -0.002 &-0.058  &0.964
                             \end{bmatrix}\]
The autoregressive model of $\YS_2(t)$ can be considered as a \textit{robust} stock return volatility model - LPVolatility Modeling, which is less affected by unusually large extreme events. The model for $\YS_2(t)$ automatically discover many known facts: (a) the sign of the coefficient linking volatility and return is negative -confirming ``leverage effect''; (b) $\YS_2(t)$ is positively autocorrelated, known as volatility clustering; (c) The positive interaction with lagged $ \YS_4(t)$ accounts for the ``excess kurtosis''.


\section{Conclusion} 
This article provides a pragmatic and comprehensive framework for nonlinear time series modeling that is easier to use, more versatile and has a strong theoretical foundation based on recently developed theory on unified algorithms of data science via LP modeling \citep{D12b,D13a,D13USA,Deep14LP,deep16LSSD,Deep17LPMode}. Summary and broader implications of the proposed research:
\vskip.15em
$\bullet$ From the theoretical standpoint, the unique aspect of our proposal lies in its ability to \emph{simultaneously} embrace and employ spectral domain, time domain, quantile domain, and the information domain analysis for enhanced insights, which to the best of our knowledge \emph{has not appeared in the nonlinear time series literature before}.
\vskip.15em

$\bullet$  From a practical angle, the novelty of our technique is that it permits us to use the techniques from linear Gaussian time series to create non-Gaussian nonlinear time series models with highly interpretable parameters. This aspect makes \texttt{LPTime} \emph{computationally extremely attractive} for the data scientists as they can now borrow all the standard time series analysis machinery from \texttt{R} libraries for implementation purpose.
\vskip.15em

$\bullet$ From the pedagogical side, we believe that these concepts and methods can easily be augmented with the standard time series analysis course to \emph{modernize the current curriculum} so that students can handle complex time series modeling problems \citep{em10} using the tools they are already familiar with.

\vskip.25em
The main thrust of this article is to describe and interpret the steps of \texttt{LPTime} technology to create a realistic general-purpose algorithm for empirical time series modeling. In addition, many new theoretical results and diagnostic measures are presented which laid the foundation for the algorithmic implementation of \texttt{LPTime}. We showed how LPTime can systematically \emph{explore the data} to discover empirical facts hidden in time series. For example, LPTime empirical modeling of S\&P 500 return data reproduces the `stylized facts': (a) heavy tails; (b) non-Gaussian;(c) nonlinear serial dependence; (d) tail-correlation; (e) asymmetric dependence; (f) volatility clustering; (g) long-memory volatility structure; (h) efficient market hypothesis; (i) leverage effect; (j) excess kurtosis,  in a coherent manner \emph{under single general unified framework}. We have emphasized how the statistical parameters of our model can be interpreted in the light of established economic theory.
\vskip.25em
We have recently applied this theory for large-scale eye-movement pattern discovery problem, which came out as the winner (among {\bf 82} competing algorithms) of the 2014 IEEE International Biometric Eye Movements Verification and Identification Competition \citep{deep16LPiTrack}.
\vskip.25em

We conclude with some general references. Few popular articles \cite{granger03,granger93,grangerlin,engle82,parzen79,parzen67,tukey80ts,brillinger77, brillinger04,wall12}; Books \cite{terasvirta2010book,woodward2011book,tsaybook}; and Review articles \cite{granger98NSF,hendry11}.
\vskip.25em
The proposed algorithm is implemented in the \texttt{R} package \texttt{LPTime}, which is available on \texttt{CRAN}.
\addcontentsline{toc}{section}{5~\, References}
\bib

\end{document}